\RequirePackage{fix-cm}
\documentclass[smallextended]{svjour3}       
\smartqed  
\usepackage{graphicx}
\usepackage{amsmath}        
\usepackage{amssymb}
\usepackage{mathrsfs}
\usepackage[plain]{algorithm}
\usepackage{algpseudocode}
\usepackage{multirow}
\usepackage{float}
\usepackage{subfig}
\usepackage{caption}
\usepackage[sort&compress]{natbib}
\setcitestyle{authoryear,open={(},close={)}}

\begin{document}

\title{MDE-ITMF and DEwI: Two New Multiple Solution Algorithms for Multimodal Optimization}

\titlerunning{Multiple Solution Algorithms in Multimodal Optimization}        

\author{Vin\'{i}cius Magno de Oliveira Coelho         \and \\
        Gustavo Barbosa Libotte \and \\
        Francisco Duarte Moura Neto \and \\
        Gustavo Mendes Platt \and \\
        Fran S\'{e}rgio Lobato}

\authorrunning{Coelho et al.} 

\institute{Vin\'{i}cius Magno de Oliveira Coelho \and Gustavo Barbosa Libotte \and Francisco Duarte Moura Neto \at
              Polytechnic Institute, Rio de Janeiro State University, 25, Bonfim St., Nova Friburgo, 28625--570, Brazil \\
              \email{vcoelho@iprj.uerj.br, {gustavolibotte@iprj.uerj.br}, fmoura@iprj.uerj.br}           
           \and
           Gustavo Mendes Platt \at
           School of Chemistry and Food, Federal University of Rio Grande, Cel. Francisco Borges de Lima St., 3005, Santo Ant\^{o}nio da Patrulha, 95500--000, Brazil
           \\\email{gmplatt@furg.br}
           \and
           Fran S\'{e}rgio Lobato \at
           Chemical Engineering Faculty, Federal University of Uberl{\^{a}}ndia, João Naves de Ávila Av., 2121, 38408--144, Uberl{\^{a}}ndia, Brazil
           \\\email{fslobato@ufu.br}
}

\date{Received: date / Accepted: date}

\maketitle

\begin{abstract}
Mathematical formulations of real world optimization studies frequently present characteristics such as non-linearity, discontinuity and high complexity. This class of problems may also exhibit a high number of global minimum/maximum points, especially for optimization problems arising from nonlinear algebraic systems (where null minima correspond to the solutions of the original algebraic system). Due to the multimodal nature of these functions, multipopulation methods have been employed in order to obtain the highest number of points of global minimum/maximum. In this work, two new approaches were analyzed, employing an iterative penalization technique and a multipopulation procedure -- together with the Differential Evolution algorithm -- devoted to obtain the full set of solutions for multimodal optimization problems. The first method proposed is the Multipopulation Differential Evolution with iterative technique of modification of the objective function, MDE-ITMF, and the second method proposed is the Differential Evolution with Initialization, DEwI. In this second proposal, the MDE-ITMF method is used as an initializer of the initial populations and from a given moment the Differential Evolution is used to solve the problem at hand. In both approaches, subpopulations evolve simultaneously throughout the iterative process. MDE-ITMF and DEwI methods were applied in a set of ten multimodal benchmark functions. Based on the results obtained, we can conclude that MDE-ITMF and DEwI are suitable and promising tools for multimodal optimization.
\end{abstract}

\section{Introduction}
\label{intro}

Commonly, engineering models give rise to optimization problems with highly nonlinear functions, great complexity and a large number of solutions. Due to these characteristics, stochastic methods have been widely used to solve such problems, as presented in \cite{bib:alvarez2008parameter, bib:bonilla2009calculation, bib:platt2013calculation, bib:platt2014computational}. Among the existing stochastic methods, Differential Evolution (DE), which was proposed by \cite{bib:price1995differential}, has been widely used in the last years. The works of \cite{bib:ponsich2011differential, bib:li2012application, bib:dragoi2016use} are examples of the proven effectiveness of DE in problems with relative complexity.

Researchers have continuously proposed changes to stochastic algorithms in order to increase the effectiveness of such techniques, as well as to reduce the computational cost. In this context, \cite{bib:zaharie2004multipopulation} proposed the Multi-resolution Multipopulation Differential Evolution (MMDE), which employs a subpopulation and a multi-resolution approach, aiming to ensure the achievement of global optima. For each subpopulation, standard genetic operators of DE are applied in order to increase the convergence speed of the algorithm. MMDE uses a multi-resolution approach to ensure that different optimum points are achieved. Nevertheless, there are no restrictions on the exploitation of subpopulations, which means that different subpopulations may converge to the same point. \cite{bib:thomsen2004multimodal} proposed Crowding Differential Evolution (CDE) and Sharing Differential Evolution (ShDE) techniques to solve multimodal optimization problems. Basically, CDE replaces an offspring with the most similar individual among the individuals in the subset formed by an agglomeration parameter, while ShDE uses a function to share the value of the objective function between similar individuals. The main objective of such procedures is to group similar individuals, in order to obtain all solutions from multimodal optimization problems. \cite{bib:thomsen2004multimodal} pointed out that CDE obtained better results in comparison to DE and ShDE, as it obtained a higher number of optimum points and required a smaller number of evaluations of the objective function to satisfy the same stopping criterion. Species-based Differential Evolution (SDE), proposed by  \cite{bib:li2005efficient}, splits the population into several subpopulations and, at each iteration, DE operators are employed in each subpopulation, which are generated by defining an agglomeration radius around the fittest individuals. In this way, SDE is able to obtain multiple optimizers from a single run of the algorithm. \cite{bib:yu2011multi} use a multi-population approach for the generation of adaptive parameters, aiming to increase the diversity of the members of the subpopulations in the evolutionary method. \cite{bib:ali2015multi} proposed the Multi-Population Differential Evolution with balanced ensemble of mutation strategies (mDE-bES), where each subpopulation uses different mutation strategies to increase the diversity and share information with each other. In turn, \cite{bib:cui2016adaptive} proposed the Adaptive Differential Evolution using Multiple sub-populations (MPADE), where the population is subdivided in terms of the value of the objective function and the authors make use of adaptive parameters in an attempt to obtain several optimums simultaneously. More recently, \cite{bib:Liao2018} proposed the Dynamic Repulsion-based Evolutionary Algorithms (DREA) that integrates dynamic repulsion techniques and evolutionary algorithms. DREA uses an adaptive parameter that controls the repulsion radius in each external iteration. Essentially, the difference between DREA and the methods proposed here is mainly due to the fact that the evolutionary strategies proposed in this work cause subpopulations to evolve simultaneously, whereas DREA generates a new initial population at the beginning of each external iteration.

The objective of this work is to propose two evolutionary techniques, based on DE, for calculating optimal solutions of multimodal optimization problems: Multipopulation Differential Evolution with Iterative Technique of Modification of the Objective Function (MDE-ITMF) and Differential Evolution with Initialization (DEwI). Both methods are analyzed taking into account a set of benchmark functions and the obtained results are compared to those calculated sequentially by the canonical version of DE. We also carried out a study on the control parameters of each method, in order to analyze the evolutionary behavior for different sets of parameters. Information related to the execution time, number of evaluations of the objective function and number of distinct global optimum points obtained are analyzed. Based on these values, a statistical study on the result is performed. It is observed a considerable gain of both methods proposed in this work in relation to DE (when executed sequentially), especially regarding the real behavior of the simultaneous evolution of the subpopulations.

The remainder of this work is structured as follows: DE is briefly described in Section~\ref{sec:DE}. Sections~\ref{sec:MDEITMF} and \ref{sec:DEwI} are intended for the proposals of the Multipopulation Differential Evolution with Iterative Technique of Modification of the Objective Function and the Differential Evolution with Initialization methods, respectively, which represent the major contributions of this work. Section~\ref{sec:results} presents the numerical results. Finally, conclusions are outlined in Section~\ref{sec:conclusions}.

\section{Description of the algorithms}
\label{sec:novel_multmod_method}

In this section, a brief description of the Differential Evolution algorithm is presented. We also describe two novel multimodal and multipopulation algorithms. The first proposal is the method of Multipopulation Differential Evolution with Iterative Technique of Modification of the Objective Function (MDE-ITMF). The second proposal uses the developed method, MDE-ITMF, as an initialization tool for the initial populations. Thus, after a certain point, the canonical Differential Evolution method is used to obtain the global minimum point of the approached function. This second approach is called Differential Evolution with Initialization method (DEwI).

\subsection{Differential Evolution -- DE}
\label{sec:DE}

The Differential Evolution (DE) is a direct search method that uses a population of individuals in its iterative processes in such a way that the most able individuals succeed for the next iterations. This method has been widely used in the field of optimization with numerous applications in engineering problems; see, for instance \cite{bib:4580136}, \cite{bib:4441675}, \cite{bib:lobato2008}, \cite{bib:SACCO20091093}, \cite{bib:Regulwar2010} and \cite{bib:LOBATO2010}.

We will denote by $\boldsymbol{x}$ a generic element in the domain of the objective function; typically, $\boldsymbol{x} \in \mathbb{R}^{d}$, where $d$ represents the dimension of the problem. The population evolution proposed by DE follows three fundamental steps: mutation, crossover and selection. The optimization process starts by creating a collection of individuals (points in $\mathbb{R}^{d}$) called the initial
3
 population, consisting of $N_{\text{p}}$ elements, represented by the set

\begin{equation}
\label{eq:pop_inicial}
X^{(0)} = \{\boldsymbol{x}_{i} ^{(0)} \in \mathbb{R}^{d}, \; i = 1, \dots, N_{\text{p}}\} \: .
\end{equation}

The individuals in the initial population are chosen randomly in the function domain. In the first step, the mutation operator creates a mutant vector, adding the balanced difference between two individuals to a third member of the current population, as follows,

\begin{equation} \label{eq:mutacao_DE_canonica}
\boldsymbol{v}_{i}^{\left(G+1\right)} = \boldsymbol{x}_{r_{1}}^{\left(G\right)} + F\left(\boldsymbol{x}_{r_{2}}^{\left(G\right)} - \boldsymbol{x}_{r_{3}}^{\left(G\right)}\right) \; ,
\end{equation}

\noindent where $i = 1, \dots, \; N_{\text{p}}$. The three distinct indices, denoted by $ r_{1} $, $ r_{2} $ and $ r_{3} $, are randomly selected in the individuals of the population $N_{\text{p}}$. The parameter $ F $ represents the amplification factor, which controls the contribution added by the vector difference. According to \cite{bib:storn1997differential}, $ F \in \left[0,\;2\right] $. In this work, the value of $F$ was contained in the interval $[0,1]$.

The second step of the algorithm is the application of the crossover process, intended to increase the diversity of the trial vector ($\boldsymbol{u}_{i}$), by modifying its current coordinates. This genetic operator creates new candidates by combining the attributes of the individuals of the original population with those resulting in the mutation step, through the operation

\begin{equation} \label{eq:crossover_DE_canonica}
 u_{i} ^{\left(G+1\right)} (k) =
 \begin{cases}
 v_{i} ^{\left(G+1\right)} (k), \; \text{if} \; \text{rand}\left(k\right) \leq \text{\emph{CR}} \; \text{or} \; k = \text{rnbr}\left( i \right)\\
 x_{i} ^{\left(G\right)} (k), \; \text{if} \; \text{rand}\left(k\right) > \text{\emph{CR}} \; \text{and} \; k \neq \text{rnbr}\left( i \right)
 \end{cases} \; ,
 \end{equation}
 
\noindent where $x_ {i} (k)$ denotes the \textit{k}-th entry of the vector $ \boldsymbol{x}_ {i}$, $k = 1, \dots, \; d$ and $ \text{rand} \left(k \right) \in \left[0, \; 1 \right] $ is a random real number with uniform distribution. The choice of the attributes of a given individual is defined by the crossover coefficient, represented by $ CR $, such that $ CR \in \left[0,\;1\right] $ is a parameter defined by the user. The term $rnbr \left( j \right) \in \left[1,\;d\right]$ is a randomly chosen index and, according to \cite{bib:storn1997differential}, its purpose is to ensure that at least one element of $\boldsymbol{v}_{i}^{\left(G+1\right)}$ will be selected to form the trial vector $\boldsymbol{u}_{i}^{\left(G+1\right)}$.

The selection step ensures that the best individuals succeed for the next generation. For this purpose, at the end of each iteration (generation) a comparison is made between the trial vectors, represented by Eqs.~\ref{eq:crossover_DE_canonica}, and the individuals of the current population. This comparison is made from the value of the objective function of each individual. The individuals who obtain better values of the objective function, in that comparison, form the population that will go on to the next generation, thus continuing the iterative process of the algorithm.

Explicitly, if ${\cal {F}}(\boldsymbol{u}_{i}^{(G + 1)}) \leq {\cal{F}}(\boldsymbol{x}_{i}^{(G)})$, then $\boldsymbol{x}_{i}^{(G + 1)} = \boldsymbol{u}_{i}^{(G + 1)}$; otherwise, the individual of the population remains in the next iteration, in other words, $\boldsymbol{x}_{i}^{(G + 1)} = \boldsymbol{x}_{i}^{(G)}$, where ${\cal {F}}$ denotes the objective function.

The typical stopping criterion for the Differential Evolution algorithm is to consider a maximum number of iterations, $G_{\text{max}}$. In addition, a measure of spreading of the population of the algorithm will be used in this work, denoted by ${\cal{D}}^{(G)}$, and is stipulated to be less than a given quantity $\epsilon > 0$.

Thus, the algorithm is terminated if either of the two criteria is met. To define the spreading measure, the element of the population that minimizes the objective function is chosen in the current generation,

\begin{equation}
\boldsymbol{x}^{(G)} = \underset{i=1, \dots, N_{\text{p}}}{\text{arg min}} \; {\cal{F}}(\boldsymbol{x_{i} ^{(G)}}) \: ,
\end{equation}

\noindent and the average relative normalized distance between all individuals of the $ G $ generation and the individual of the population that minimizes the objective function is calculated,

\begin{subequations}
\begin{equation}
\label{eq:dist_norm_rel_ED_gener}
{\cal{D}}^{(G)} = \frac{1}{N_{\text{p}}}  \sum_{i=1}^{N_{\text{p}}} \frac{\left( \sum_{k=1}^{d} \left( \frac{x_{i} ^{(G)} (k) - x ^{(G)} (k)}{U_{k} - L_{k}} \right)^{2} \right)^{\frac{1}{2}} } { \left( \sum_{k=1}^{d} \left( \frac{x ^{(G)} (k)}{U_{k} - L_{k}} \right)^{2} \right)^{\frac{1}{2}} } \; ,
\end{equation}

\noindent where $L = (L_{1}, \cdots, L_{d})$ and $ U = (U_{1}, \cdots, U_{d})$ represent the upper and lower limits of the parallelogram representing the function domain. When $d = 2$, for instance, ${\cal{D}}^{(G)}$ is reduced to

\begin{equation}
\label{eq:dist_norm_rel_ED}
{\cal{D}}^{(G)} = \frac{1}{N_{\text{p}}}  \sum_{i=1}^{N_{\text{p}}} \frac{{\sqrt{ \left( \frac{x_{i} ^{(G)} (1) - x ^{(G)}(1)}{L_{x}} \right)^{2} + \left( \frac{x_{i} ^{(G)} (2) - x ^{(G)} (2)}{L_{y}} \right)^{2} }}} {\sqrt{ \left( \frac{x ^{(G)} (1)}{L_{x}} \right)^{2} + \left( \frac{x^{(G)} (2)}{L_{y}} \right)^{2} }} \; ,
\end{equation}
\end{subequations}

\noindent where the parameters $L_{x}$ and $L_{y}$ are the sides, respectively, in the abscissa and ordinate, of a rectangle that contains the function domain, and where, for example, the vector $(x_{i} ^{(G)} (1), x_{i} ^{(G)} (2))$ denotes the $i$-th individual of the $G$-th iteration. Particularly, $L_{x} = U_{1} - L_{1} \; \text{and} \; L_{y} = U_{2} - L_{2} $.

\subsection{Multipopulation Differential Evolution with Iterative Technique of Modification of the Objective Function -- MDE-ITMF}
\label{sec:MDEITMF}

The main objective of the algorithm proposed here is to provide a new multipopulation framework for the Differential Evolution method, in order to create a new procedure that is capable of obtaining, simultaneously, all the global solutions of an optimization problem in a single execution of the algorithm.

The main part of the methodology is called the \textit{iterative modification of the objective function} --- which effectively constitutes a penalty added to the objective function --- and, together, the multipopulation approach is used so that the algorithm is able to calculate all global solutions of the optimization problem simultaneously, considering that each subpopulation distinguishes a different objective function due to the application of the penalty/modification technique. Furthermore, there is a natural repulsion between the subpopulations, as a consequence of the penalty functions. As an illustration, in a problem with four subpopulations, each population distinguishes an objective function with three different ``repulsion areas'', as a consequence of the penalties for the other three subpopulations.

Through this simple approach, the diversity of the population, as a whole,  is preserved, which can increase computational gain and accelerate convergence. In addition, the algorithm naturally becomes parallelizable, which guarantees the possibility of increased performance.

Given an optimization problem

\begin{equation}
\underset{\boldsymbol{x} \in D}{\min} \; {\cal{F}}(\boldsymbol{x}),
\end{equation}

\noindent where ${\cal{F}}: \mathbb{R}^{d} \rightarrow \mathbb{R}$ is the objective function to be optimized, and assume that the feasible domain is given by

\begin{equation} 
D = \{\boldsymbol{x} \in \mathbb{R}^{d} \mid g_{l} (\boldsymbol{x}) \leq b_{l} (\boldsymbol{x}), \; \text{with} \; l = {1, \dots, r}\} \: . 
\label{eq:dxe}
\end{equation}

In this sense, $\boldsymbol{x} = (x(1), \dots, x(d))^{\top}$ is the optimization variable.

The process of obtaining the optimal global solutions begins with the creation of all individuals in the initial population, which is contained in the feasible domain $D$. Initially, the number of subpopulations $N_{\text{sp}}$ that will compose the final population, besides the number of individuals $N_{\text{p}}$ that integrate each subpopulation, are defined. The idea is that the population consists of $N_{\text{p}} \times N_{\text{sp}}$ individuals, where each subpopulation will be managed independently of the others. Thus, at the end of $ G_{\text{max}} $ generations, each of the $N_{\text{p}}$ subpopulations is expected to converge to a distinct global optimal solution.

Denote by

\begin{equation}
\boldsymbol{x}_{ij} ^{(G)} \in \mathbb{R}^{d} \: ,
\end{equation}

\noindent the \textit{i}-th indivídual, $i = 1, \dots, N_{\text{p}}$, of \textit{j}-th subpopulation, $j = 1, \dots, N_{\text{sp}}$, on \textit{G}-th iteration of algorithm. The entries of vector $\boldsymbol{x}_{ij} ^{(G)}$ are denoted by

\begin{equation}
x_{ij} ^{(G)} (k), \; k = 1, \dots, d \:.
\end{equation}

\noindent where $d$ represents the dimention of the problem.

Thus, a three-dimensional matrix, $\mathbb{X}$, is constructed to represent the subpopulations,

\begin{equation}
(\mathbb{X}^{(G)})_{kij} = x_{ij} ^{(G)}(k) \: , 
\label{eq:Xtridimensional}
\end{equation}

\noindent where $k = 1, \dots, d, \: i = 1, \dots, N_{\text{p}}, \: \text{e} \: j = 1, \dots, N_{\text{sp}}$.

\begin{figure}[!htb]
\includegraphics[scale=0.5]{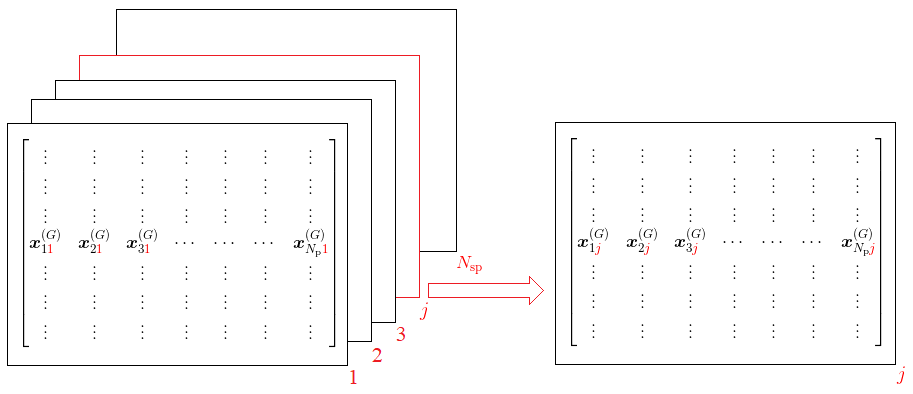}
\caption{Structuring the three-dimensional matrix $ \mathbb{X} $ according to subpopulations}
\label{fig:matrizTridimensional_Xestruturada}
\end{figure}

The Figure~\ref{fig:matrizTridimensional_Xestruturada} presents the structure of the matrix $\mathbb{X}$ according to subpopulations. On the right side of this figure, for example, $ \boldsymbol{x}^{(G)} _{2 j} $, the second column of the indicated matrix, represents the second individual in the $ j $-th subpopulation in $G$-th generation.

The initialization of individuals is performed randomly in $D$. Therefore, search intervals with lower $(\boldsymbol{L})$ and upper $(\boldsymbol{U})$ bounds are pre-established in each of the $ d $ dimensions of the problem, thus generating a hyperspace in the form of a hyper-parallelepiped in $\mathbb{R}^{d}$, given by the constraints

\begin{equation}
L_{k} \leq x(k) \leq U_{k}, \; \text{for all} \; k = 1, \dots, d,
\end{equation}

\noindent with $\boldsymbol{x} = (x(1), \dots ,x(d) )^{\top} \in \mathbb{R}^{d}$, where $L_{k}$ and $U_{k}$ are pre-established. 

In this way, the coordinates of the population elements can be represented by

\begin{equation} \label{eq:init_population}
x_{ij}^{\left(0\right)} (k) = L_{k} + h \left( U_{k} - L_{k} \right) \; ,
\end{equation}

\noindent with $i = 1, \dots, \; N_{\text{p}}$, $j = 1, \dots, \; N_{\text{sp}}$ and $k = 1, \dots, \; d$. To ensure that each individual is within the expected search domain, $h$ is a random number with uniform distribution in the interval $\left[ 0, \; 1 \right]$.

The crucial procedure that allows the simultaneous obtaining of distinct optimal points occurs through the creation of areas of repulsion around the best approximations of the current generation. Among the many techniques available in the literature for this same purpose, the deflation strategy proposed in \cite{bib:hirsch2009solving} was used. Due to the fact that repeated executions of a stochastic optimization algorithm do not guarantee that different solutions will be obtained in a multimodal optimization problem, this technique is designed to search for approximations that were not previously obtained by other subpopulation. It should be noted that all $ N_{\text{sp}} $ subpopulations evolve concurrently over the generations.

Initially, the best approximations of each subpopulation are selected. Be $ S^{(G)} $ a matrix $ d \times N_{\text{sp}} $ whose columns are the best approximations of the solutions to the problem, denoted by $ \boldsymbol{s}_{ _{1}} ^{(G)}, \dots, \boldsymbol{s}_{N_{\text{sp}}} ^{(G)} $, obtained by each subpopulation, in the generation $G$. Thus, 

\begin{equation}
 S^{(G)} = (\boldsymbol{s}_{ _{1}} ^{(G)}, \dots, \boldsymbol{s}_{ _{N_{\text{sp}}}} ^{(G)}) \in \mathbb{R}^{d}.
\end{equation}

In this way,

\begin{equation}
\boldsymbol{s}_{j} ^{(G)}= \underset{\boldsymbol{x}^{(G)} _{\tilde{\imath}j}, \: \tilde{\imath} = 1, \dots, N_{\text{p}}}{\text{arg min}} \; {\cal{F}}\left(\boldsymbol{x}_{\tilde{\imath}j}^{\left(G\right)}\right) \; ,
\label{eq:formulacao_S}
\end{equation}

\noindent where $j = 1, \dots, N_{\text{sp}}$.

Let $j$ be the number of individuals in the matrix $S$, thus $j = 1, \dots, N_{\text{sp}}$. The penalty technique creates a repulsion area in the objective function around the approximate solutions, except the solution $\boldsymbol{s}_{j}$, considering the modification

\begin{equation} \label{eq:funcao_penalizacao}
{\cal{F}}_{j, \beta, \rho}\left(\boldsymbol{x}\right) = {\cal{F}}(\boldsymbol{x}) + \beta\sum\limits_{\kappa=1, \kappa\neq j}^{N_{\text{sp}}}\text{e}^{-\Vert \boldsymbol{x} - \boldsymbol{s}_{\kappa} \Vert}\Phi_{\rho}\left(\Vert \boldsymbol{x} - \boldsymbol{s}_{\kappa} \Vert\right) \: ,
\end{equation}

The constants $ \beta $ and $ \rho $ represent the magnitude of the penalty and the limit of the penalization region, respectively. The function represented by $\Phi_{\rho}(\delta)$ is employed to activate/deactivate the penalty technique, where $\delta$ represents the distance between a generic individual $\boldsymbol{x}$ and the best minimizer previously obtained
$\boldsymbol{s}_{\kappa}$. In this sense, the  $ \Phi_{\rho}\left(\delta\right) $ function returns only two values, as shown below

\begin{equation}
\Phi_{\rho}(\delta) = 
\begin{cases}
1, \; \text{if} \; \delta \leq \rho \\
0, \; \text{otherwise}
\end{cases} \: .
\end{equation}

Thus, $ \beta $ must be a value that generates an appropriate penalty in relation to the value of the objective function of the problem, while $ \rho $ must be a small value in such a way that it does not generate a penalty area where there are solutions that have not yet been obtained.

In order to use the penalty technique in a system of non-linear equations, it is necessary to transform the system into an appropriate (and scalar) objective function, by summing the squares of the residues of each nonlinear equation, for example, as follows, 

\begin{equation} \label{eq:modelo_fitness}
{\cal{F}}\left(x\right) = \sum\limits_{\alpha = 1}^{n} f_{\alpha}^{2}\left(x\right) \; ,
\end{equation}

\noindent where $ n $ is the number of equations of the system. Thus, a possible function that represents the penalty technique applied to a system of nonlinear equations can be represented as

\begin{equation} \label{eq:funcao_penaliz_sist}
{\cal{F}}_{j, \beta, \rho}\left(\boldsymbol{x}\right) = \sum\limits_{\alpha = 1}^{n} f_{\alpha}^{2}\left(x\right) + \beta\sum\limits_{\kappa=1, \kappa\neq j}^{N_{\text{sp}}}\text{e}^{-\Vert \boldsymbol{x} - \boldsymbol{s}_{\kappa} \Vert}\Phi_{\rho}\left(\Vert \boldsymbol{x} - \boldsymbol{s}_{\kappa} \Vert\right) \: .
\end{equation}

Algorithm~\ref{alg:penaliz} presents a possible implementation of the penalty function described above.

\begin{algorithm}[!htb]
\caption{The penalty technique applied to an objective function}\label{alg:penaliz}
\begin{algorithmic}[1]
\Function{modFun}{$\boldsymbol{x},\; j, \; S , \; \beta, \; \rho$}
\State $ \Lambda \leftarrow 0 $
\For{\textbf{each} $\kappa \neq j$} 
\State $ \delta \leftarrow \Vert \boldsymbol{x} - \boldsymbol{s}_{\kappa} ^{(G)} \Vert $
\State Calculate $ \Phi_{\rho}\left(\delta\right) $
\State $ \Lambda \leftarrow \Lambda + \text{e}^{-\delta}\Phi_{\rho}\left(\delta\right) $
\EndFor
\State $ y \leftarrow {\cal{F}}\left(\boldsymbol{x}\right) + \beta\Lambda $\Comment{According to Equation~\ref{eq:funcao_penalizacao}}
\State \Return $ y $
\EndFunction
\end{algorithmic}
\end{algorithm}

Essentially, the proposed method, MDE-ITMF, make use of the three genetic operators from DE, as described in Section~\ref{sec:DE}, for each subpopulation. The main difference, however, is in the Selection operator. The objective function is penalized, in order to guarantee that all optimum points are obtained, by each subpopulation. Such procedure guarantees the obtaining of multiple optima in a single execution.

Mutation and recombination steps are similar to those used in the canonical version of DE; however, there is a condition so that the individuals generated by these procedures must be in the domain of the problem.

Algorithm~\ref{alg:etapa_selecao_MDEITMF} presents the selection procedure of MDE-ITMF.

\begin{algorithm}[!htb]
\caption{Selection step with iterative modification of the objective function}\label{alg:etapa_selecao_MDEITMF}
\begin{algorithmic}[1]
\Function{selectionStep}{$ \boldsymbol{x}_{c}, \;\boldsymbol{u}, \;j , \;N_{\text{p}}, \;d,\;S, \;\beta, \;\rho, \;\boldsymbol{L}, \;\boldsymbol{U} $}
\For{$ i \leftarrow 1 : N_{\text{p}} $}
\State $ E \leftarrow \textbf{true} $\Comment{Flag to control individual evolution}
\For{$ k \leftarrow 1 : d $}
\If{$ \boldsymbol{u}_{i}(k) < L_{k} $ \textbf{or} $ \boldsymbol{u}_{i}(k) > U_{k} $}
\State $ E \leftarrow \textbf{false} $
\State break
\EndIf
\EndFor
\If{$ E = \textbf{true} $}
\If{\Call{modFun}{$\boldsymbol{u}_{i}, \;j, \;S, \;\beta, \;\rho$} $ < $ \Call{modFun}{$\boldsymbol{x}_{ci}, \;j ,\;S,\;\beta,\;\rho$}}
\State $ \boldsymbol{x}_{ci} \leftarrow \boldsymbol{u}_{i} $
\EndIf
\EndIf
\EndFor
\State \Return $ \boldsymbol{x}_{c} $
\EndFunction
\end{algorithmic}
\end{algorithm}

Using this strategy, the current subpopulation evolves considering only one specific problem, introduced by the dynamic modification of the objective function, since the matrix $S$ changes at each generation. Thus, the algorithm is oriented in a way so that two subpopulations will not converge to the same approximation of the solution, and can evolve independently and simultaneously.

\begin{figure}[!htb]
\centering
\subfloat[][$ 1^{\mathrm{st}} $ subpopulation]{\label{fig:G1P1}
\includegraphics[scale=0.41]{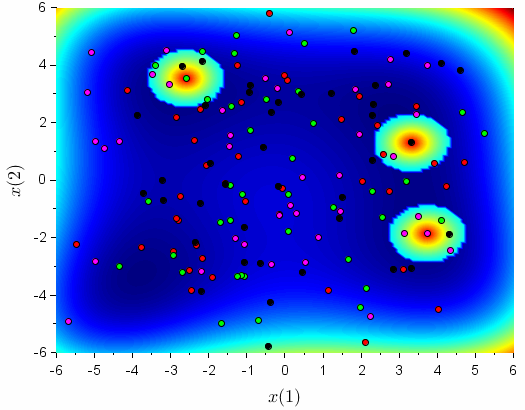}} \;
\subfloat[][$ 2^{\mathrm{nd}} $ subpopulation]{\label{fig:G1P2}
\includegraphics[scale=0.41]{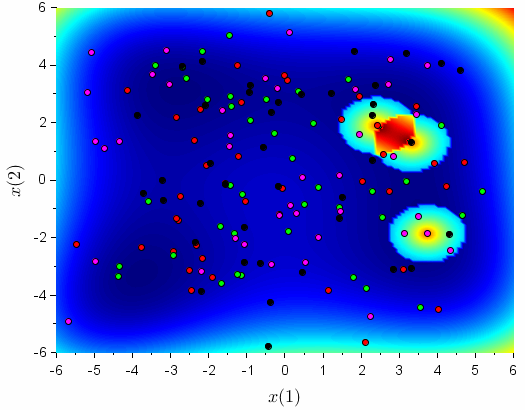}} \\  
\subfloat[][$ 3^{\mathrm{rd}} $ subpopulation]{\label{fig:G1P3}
\includegraphics[scale=0.41]{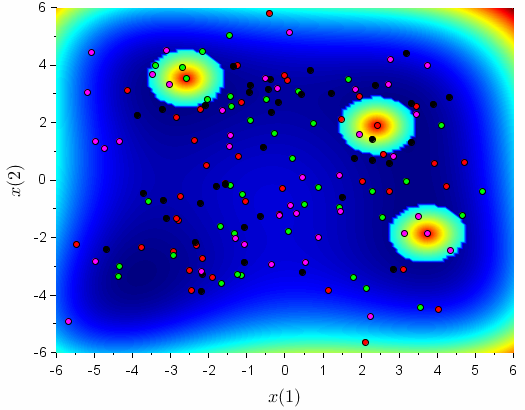}} \;
\subfloat[][$ 4^{\mathrm{th}} $ subpopulation]{\label{fig:G1P4}
\includegraphics[scale=0.41]{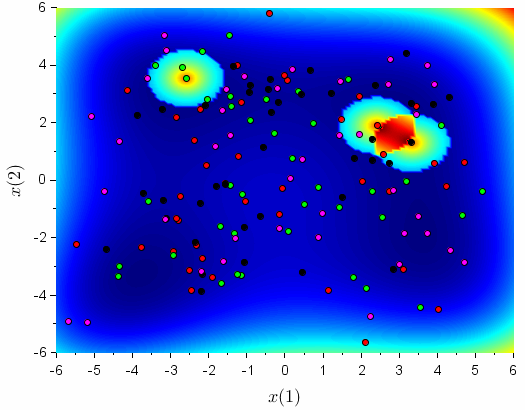}}
\caption{Evolution of each subpopulation for $ G = 1 $}
\label{fig:exemplo1}
\end{figure}

Figure~\ref{fig:exemplo1} contains a pictorial representation of the modification of the objective function, using the MDE-ITMF method in the Himmelblau function during the first iteration, in relation to each subpopulation. The subpopulations are marked by different colors (red, green, black and pink). This figure was generated using the following set of parameters: $N_{\text{sp}} = 4$, $N_{\text{p}} = 40$, $F = 0,4$, $\emph{CR} = 0,1$ $\rho = 1$ and $\beta = 2 \times 10^{3}$. The values of parameters $F$ and \emph{CR} were chosen taking into account the suggestion presented by \cite{bib:storn1997differential}.

The values for the parameters $N_{\text{sp}}$, and $N_{\text{p}}$ were chosen, respectively, from the number of global minimum points in the Himmelblau function (four minima) and the number of individuals that will be illustrated (for each subpopulation). With this set of parameters, one can observe a ``picture'' for each subpopulation in the initial generation of the process. Clearly, the individuals of the subpopulations are randomly distributed in the feasible domain of the problem. On the other hand, we can see different objective functions for the subpopulations, which indicates that these subpopulations will mutually repel in the subsequent generations.

\begin{figure}[!htb]
\centering
\subfloat[][$ 1^{\mathrm{st}} $ subpopulation]{\label{fig:G27P1}
\includegraphics[scale=0.41]{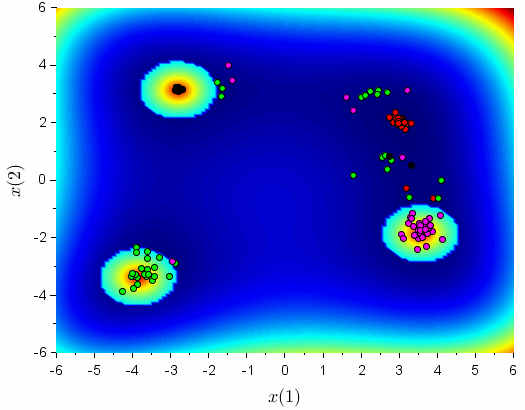}} \;
\subfloat[][$ 2^{\mathrm{nd}} $ subpopulation]{\label{fig:G27P2}
\includegraphics[scale=0.41]{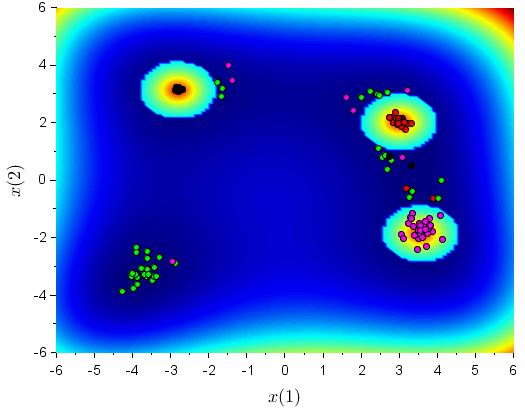}} \\  
\subfloat[][$ 3^{\mathrm{rd}} $ subpopulation]{\label{fig:G27P3}
\includegraphics[scale=0.41]{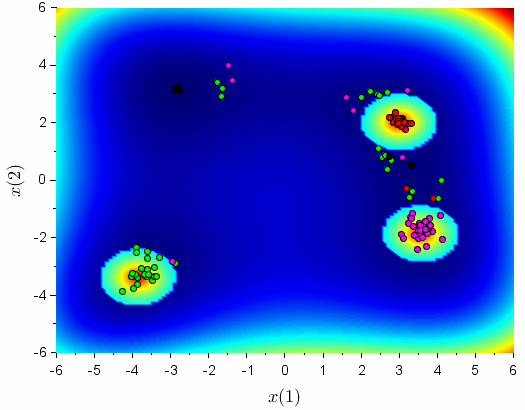}} \;
\subfloat[][$ 4^{\mathrm{th}} $ subpopulation]{\label{fig:G27P4}
\includegraphics[scale=0.41]{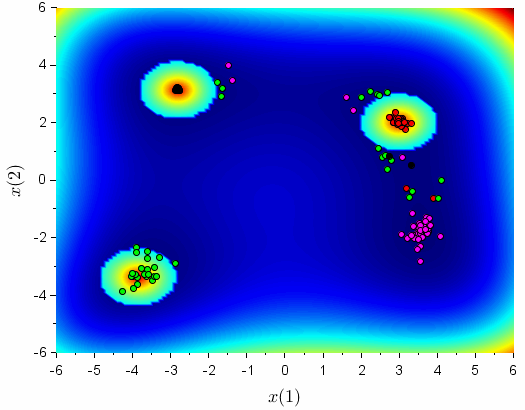}}
\caption{Evolution of each subpopulation for $ G = 27 $}
\label{fig:exemplo2}
\end{figure}

Figure~\ref{fig:exemplo2} represents a different stage in the evolutionary process of the MDE-ITMF method presented in Figure~\ref{fig:exemplo1}, now in its twenty-seventh iteration. This figure clearly indicates the convergence areas of the different subpopulations taking into account the different points of global optimum. It is important to note that during the iterations the individuals of the subpopulations converged to the same point of global optimum, and in the course of the iterations the diversity of the subpopulations tends to decrease taking into account that the individuals do not migrate between the subpopulations. Another interesting feature that can be noted in a comparison between Figures~\ref{fig:exemplo1} and ~\ref{fig:exemplo2} is that the penalty regions are dynamical, i.e., these regions move from one generation to another. 

The stopping criterion of the algorithm obeys two distinct rules: one regarding the number of iterations of the algorithm, and another regarding the spreading measure of the subpopulations. The method ends its execution immediately when the generation counter reaches the maximum limit established a priori by the user, represented by $G_{\text{max}}$. Due to the fact that the subpopulations are independent, each of them is subject to a maximum of $ G_{\text{max}} $ generations. In addition, there is an additional stopping criterion responsible for controlling the convergence of each subpopulation. This stopping criterion aims to reduce the maximum number of operations (and, consequently, the computation time) to obtain the solutions. The stopping criterion uses the formulation presented by Eq.~\ref{eq:dist_norm_rel_ED_gener} extended to all subpopulations. In this way, the average relative normalized distance between all individuals of each subpopulation and the individual that represents the best approximation obtained by that subpopulation is calculated. If
 $d = 2$, this quantity for each subpopulation is calculated by
 
 \begin{equation} \label{eq:Dist_norm_relat_subpop}
{\cal{D}}_{j}^{(G)} = \frac{1}{N_{\text{p}}}  \sum_{i=1}^{N_{\text{p}}} \frac{{\sqrt{ \left( \frac{x_{i j} ^{(G)} (1) - s_{j} ^{(G)}(1)}{L_{x}} \right)^{2} + \left( \frac{x_{i j} ^{(G)} (2) - s_{j} ^{(G)} (2)}{L_{y}} \right)^{2} }}} {\sqrt{ \left( \frac{s_{j} ^{(G)} (1)}{L_{x}} \right)^{2} + \left( \frac{s_{j}^{(G)} (2)}{L_{y}} \right)^{2} }} \; ,
\end{equation}

\noindent where $ j = 1, \dots, \; N_{\text{sp}} $, where remember $\boldsymbol{s}_{j} ^{(G)}$ is the best individual in the subpopulation $j$. The parameters $L_{x}$ and $L_{y}$ are the sides of a rectangle that contains the function domain.

Thus, the algorithm is terminated for a $j$-th subpopulation in the $G$-th iteration if convergence is achieved, in the sense that

\begin{equation} \label{eq:dist_norm_rel_MDEITMF_epsilon}
{\cal{D}}_{j}^{\left(G\right)} < \epsilon \: .
\end{equation}

It is observed that the lower the $\epsilon$ value, with $\epsilon > 0$, the closer the individuals of the subpopulation will be to each other, offering a more refined solution.

The adoption of both stopping criteria causes a positive impact in the convergence process, since it will not be necessary to apply the entire iterative process previously described in the subpopulations that have already converged.

Consequently, the algorithm is terminated if

\begin{subequations}
\begin{equation}
\label{eq:critp_MDEITMF}
G = G_{\text{max}},
\end{equation}

\noindent or, 

\begin{equation}
{\cal{D}}_{j}^{\left(G\right)} < \epsilon \: .
\end{equation}
\end{subequations}

The main computational routine of MDE-ITMF is described in detail in the Algorithm~\ref{alg:funcao_principal_MDEITMF}.

\begin{algorithm}[!htb]
\caption{MDE-ITMF Algorithm}
\label{alg:funcao_principal_MDEITMF}
\begin{algorithmic}[1]
\Function{MDE-ITMF}{$N_{\text{sp}},\;N_{\text{p}},\;d,\;F,\;\text{\emph{CR}},\;\boldsymbol{L},\;\boldsymbol{U},\;\beta,\;\rho,\;\epsilon,\;G_{\text{max}}$}
\State $ \boldsymbol{x}^{\left(0\right)} \leftarrow $ \Call{InitPopulation}{$d,\;N_{\text{sp}},\;N_{\text{p}},\;\boldsymbol{L},\;\boldsymbol{U}$}\Comment{eq.~\ref{eq:init_population}}
\State Evaluate all individuals and get the individual with the best objective function value from each subpopulation.
\For{$ G \leftarrow 1 : G_{\text{max}} $}
\For{$ j \leftarrow 1 : N_{\text{sp}} $}
\State Compute $ {\cal{D}}_j ^{(G)} $ according to eq.~\ref{eq:Dist_norm_relat_subpop}
\If{$ {\cal{D}}_{j}^{(G)} \geq \epsilon $}
\State Initialize $ S^{(G)} \leftarrow \emptyset $\Comment{$ S^{(G)} $ contains $\left( N_{\text{sp}}\right) $ elements}
\For{$ j_{aux} \leftarrow 1 : N_{\text{sp}} $ \textbf{and} $ j_{aux} \neq j $}
\State Compute $ S^{(G)} $ according to eq.~\ref{eq:formulacao_S}

\EndFor
\State $ \boldsymbol{x}_{c}^{\left(G\right)} \leftarrow \boldsymbol{x}_{j}^{\left(G\right)} $
\State $ \boldsymbol{v} \leftarrow $ \Call{mutationStep}{$\boldsymbol{x}_{c}^{\left(G\right)},\;N_{\text{p}},\;d,\;F$}\Comment{eq.~\ref{eq:mutacao_DE_canonica}}
\State $ \boldsymbol{u} \leftarrow $ \Call{crossoverStep}{$\boldsymbol{x}_{c}^{\left(G\right)},\;\boldsymbol{v},\;N_{\text{p}},\;d,\;\text{\emph{CR}}$}\Comment{eq.~\ref{eq:crossover_DE_canonica}}
\State $ \boldsymbol{x}_{j}^{\left(G\right)} \leftarrow $ \Call{selectionStep}{$ \boldsymbol{x}_{c}^{\left(G\right)}, \;\boldsymbol{u}, \;j , \;N_{\text{p}}, \;d,\;S, \;\beta, \;\rho, \;\boldsymbol{L}, \;\boldsymbol{U} $}
\State Evaluate all individuals and get the individual with the best objective function from the current subpopulation.
\EndIf
\EndFor
\EndFor
\State \Return $ \mathbb{X} $
\EndFunction
\end{algorithmic}
\end{algorithm}

\subsection{Differential Evolution with Initialization -- DEwI}
\label{sec:DEwI}

The Differential Evolution with Initialization, denoted by the acronym DEwI, uses the MDE-ITMF to generate the initial subpopulations and then applies the canonical DE algorithm. Thus, initially, the DEwI method generates promising exploration sub-areas using the MDE-ITMF method.

Due to the fact that the MDE-ITMF method is used as a population initiator for the DE method, it became necessary to create a parameter that indicates when the MDE-ITMF method is deactivated and the DE method is activated. The average relative normalized distance ${\cal{D}}_{j} ^{(G)}$, Eq.~\ref{eq:Dist_norm_relat_subpop}, of the $j$-th subpopulation is used as the tool to disable the application of the MDE-ITMF method and start the DE algorithm in a subpopulation. A control variable represented as $tol$ is used for this purpose. The value of this parameter must be greater than the tolerance value of the spreading measure $\epsilon$. When $ {\cal{D}}_{j} ^{(G)} $ is less than the tolerance parameter, $ tol $, the procedure of the MDE-ITMF method is replaced by the standard DE method in the subpopulation that satisfies the condition.

The stopping criterion adopted in this method is similar to that adopted in the MDE-ITMF method. It should also be noted that both methods have two conditions for the stopping criterion, the first in relation to the maximum number of generations iterated by all subpopulations, $G_{\text{max}}$, and the second in relation to the tolerance of the spreading measure, $\epsilon$.

Algorithm~\ref{alg:etapa_selecao_EDI} represents the selection procedure, responsible to the switch between MDE-ITMF and DE methods. This algorithm uses the penalty technique of the objective function when performing the procedure of the MDE-ITMF method or uses the value of the objective function (without the penalty term) when executing the canonical DE method. In the Selection Step there is the \emph{Toggle} parameter that assumes two values, \emph{true} or \emph{false}. This parameter will assume the true Boolean value when the DEwI method uses the MDE-ITMF method, and the false Boolean value when the procedure is the ED method. In the latter case, this implies disabling the MDE-ITMF method procedure.

\begin{algorithm}[!h]
\caption{Selection step of DEwI method}
\label{alg:etapa_selecao_EDI}
\begin{algorithmic}[1]
\Function{selectionStep}{$ \boldsymbol{x}_{c}, \;\boldsymbol{u}, \;j , \;N_{\text{p}}, \;d,\;S, \;\beta, \;\rho, \;\boldsymbol{L}, \;\boldsymbol{U}, \;\emph{Toggle} $}
	\For{$ i \leftarrow 1 : N_{\text{p}} $}
		\State $ E \leftarrow \textbf{true} $\Comment{Flag to control individual evolution}
		\For{$ k \leftarrow 1 : d $}
			\If{$ \boldsymbol{u}_{i}(k) < L_{k} $ \textbf{or} $ \boldsymbol{u}_{i}(k) > U_{k} $}
				\State $ E \leftarrow \textbf{false} $
				\State break
			\EndIf
		\EndFor
		\If{$ E = \textbf{true} $}
			\If{$\emph{Toggle} = \textbf{true} $}
				\If{\Call{modFun}{$\boldsymbol{u}_{i}, \;j, \;S, \;\beta, \;\rho$} $ < $ \Call{modFun}{$\boldsymbol{x}_{ci}, \;j ,\;S,\;\beta,\;\rho$}}
					\State $ \boldsymbol{x}_{ci} \leftarrow \boldsymbol{u}_{i} $
				\EndIf				
			\Else
					\If{\Call{objFun}{$\boldsymbol{u}_{i}$} $ < $ \Call{objFun}{$\boldsymbol{x}_{ci}$}}
						\State $ \boldsymbol{x}_{ci} \leftarrow \boldsymbol{u}_{i} $
					\EndIf
			\EndIf
		\EndIf
	\EndFor
	\State \Return $ \boldsymbol{x}_{c} $
\EndFunction
\end{algorithmic}
\end{algorithm}

In order to represent the function that performs the penalty/modification of the objective function, the function \Call{modFun}{} is employed. This function characterizes the procedures of the MDE-ITMF method during the execution of the DEwI method. The input arguments of \Call{modFun}{} are the parameters referring to the MDE-ITMF method. However, due to the fact that the DEwI method also uses the original DE method, it became necessary to use a function that characterizes the objective function. In this way, the \Call{objFun}{} function is used in the Selection Stage during DE procedure; this function has as an input argument, the point --- individual of the population --- for which the value of the objective function is desired.

The main computational routine of DEwI is detailed in Algorithm~\ref{alg:funcao_principal}.

\begin{algorithm}[!h]
\caption{DEwI Algorithm} \label{alg:funcao_principal}
\begin{algorithmic}[1]
\Function{DEwI}{$N_{\text{sp}},\;N_{\text{p}},\;d,\;F,\;\text{\emph{CR}},\;\boldsymbol{L},\;\boldsymbol{U},\;\beta,\;\rho,\;\epsilon,\;G_{\text{max}}, \; tol$}
\State $ \boldsymbol{x}^{\left(0\right)} \leftarrow $ \Call{InitPopulation}{$d,\;N_{\text{sp}},\;N_{\text{p}},\;\boldsymbol{L},\;\boldsymbol{U}$} 
\State Evaluate all individuals and get the individual with the best objective function value from each subpopulation.
\State $ \emph{Toggle} \leftarrow \textbf{true} $
\For{$ G \leftarrow 1 : G_{\text{max}} $}
\State Compute $ {\cal{D}}_j ^{(G)} $ according to eq.~\ref{eq:Dist_norm_relat_subpop} 
\For{$ j \leftarrow 1 : N_{\text{sp}} $}
\If{$ {\cal{D}}_{j}^{(G)} \geq \epsilon $}
\If{$ {\cal{D}}_{j} ^{(G)} \geq tol $}
\State Initialize $ S^{(G)} \leftarrow \emptyset $ 
\For{$ j_{aux} \leftarrow 1 : N_{\text{sp}} $}
\If{$ j_{aux} \neq j $}
\State Calculate $ S^{(G)} $ 
\EndIf
\EndFor
\State $ \emph{Toggle} \leftarrow \textbf{true} $
\Else
\State $ \emph{Toggle} \leftarrow \textbf{false} $
\EndIf
\State $ \boldsymbol{x}_{c}^{\left(G\right)} \leftarrow \boldsymbol{x}_{j}^{\left(G\right)} $
\State $ \boldsymbol{v} \leftarrow $ \Call{mutationStep}{$\boldsymbol{x}_{c}^{\left(G\right)},\;N_{\text{p}},\;d,\;F$} 
\State $ \boldsymbol{u} \leftarrow $ \Call{crossoverStep}{$\boldsymbol{x}_{c}^{\left(G\right)},\;\boldsymbol{v},\;N_{\text{p}},\;d,\;\text{\emph{CR}}$} 
\State $ \boldsymbol{x}_{j}^{\left(G\right)} \leftarrow $ \Call{selectionStep}{$ \boldsymbol{x}_{c}^{\left(G\right)}, \;\boldsymbol{u}, \;j , \;N_{\text{p}}, \;d,\;S, \;\beta, \;\rho, \;\boldsymbol{L}, \;\boldsymbol{U}, \; \emph{Toggle} $}
\State Evaluate all individuals and get the individual with the best objective function value from the current subpopulation.
\EndIf
\EndFor
\EndFor
\State \Return $ \mathbb{X} $
\EndFunction
\end{algorithmic}
\end{algorithm}

\section{Numerical results} \label{sec:results}

The results presented in this section are related to the application of the standard DE, MDE-ITMF and DEwI methods in a set of benchmark functions usually employing in the literature. The functions used are shown in the Table~\ref{tab:bench_funct} and were obtained from \cite{bib:jamil2013literature}. More details regarding the benchmark funtions can be obtained in \cite{bib:jamil2013literature}. These functions are defined in a two-dimensional domain and are also multimodal.

\begin{table}[!htb]
\centering{}
\resizebox{\textwidth}{!}{
\begin{tabular}{cccc}
\hline
\multicolumn{2}{c}{Functions} & Domain & \# minimizers \\ \hline
$ B_{1} $ & Himmelblau & $ \boldsymbol{L} = [-6, -6]; \boldsymbol{U} = [6, 6]; $ & $ 4 $ \\
$ B_{2} $ & Trecanni & $ \boldsymbol{L} = [-5, -5]; \boldsymbol{U} = [5, 5]; $ & $ 2 $ \\
$ B_{3} $ & Six Hump--Camel & $ \boldsymbol{L} = [-3, -2]; \boldsymbol{U} = [3, 2]; $ & $ 2 $ \\
$ B_{4} $ & Cross--in--tray & $ \boldsymbol{L} = [-10, -10]; \boldsymbol{U} = [10, 10]; $ & $ 4 $ \\
$ B_{5} $ & Bird & $ \boldsymbol{L} = [-2 \pi, -2\pi]; \boldsymbol{U} = [2\pi, 2\pi]; $ & $ 2 $ \\
$ B_{6} $ & Branin RCOS & $ \boldsymbol{L} = [-5, 0]; \boldsymbol{U} = [10, 15]; $ & $ 3 $ \\
$ B_{7} $ & System of equations & $ \boldsymbol{L} = [-1, -1]; \boldsymbol{U} = [1, 1]; $ & $ 4 $ \\
$ B_{8} $ & Wayburn Seader's \#1 & $ \boldsymbol{L} = [-500, -500]; \boldsymbol{U} = [500, 500]; $ & $ 2 $ \\
$ B_{9} $ & Wayburn Seader's \#2 & $ \boldsymbol{L} = [-500, -500]; \boldsymbol{U} = [500, 500]; $ & $ 2 $ \\
$ B_{10} $ & Ackley \#3 & $ \boldsymbol{L} = [-32, -32]; \boldsymbol{U} = [32, 32]; $ & $ 2 $ \\ \hline
\end{tabular}}
\caption{Benchmark functions, domains and \# of minimizers}
\label{tab:bench_funct}
\end{table}

We employed the MDE-ITMF and DEwI methods in 100 sequential runs, thus generating 100 sets of global minimum points. The canonical DE method was performed $100 \times N_{\text{sp}}$ runs, in such a way that each $N_{\text{sp}}$ runs generated a set of global minimum points obtained by DE.

The performance measures used in this work are the elapsed time (in seconds), ET, the number of objective function evaluations, NFE, and the number of distinct global minimum points obtained, NGP. Besides these performance measures, some statistical quantities -- the arithmetic mean, $\mu$, the standard deviation, $\sigma$, and the coefficient of variation, $\upsilon$ -- are used. The coefficient of variation is a measure of dispersion that measures the variability of the data set in relative terms of the arithmetic mean. This coefficient is calculated from the ratio of the values of the arithmetic mean and the standard deviation; its representation is commonly performed in percentage terms. The lower the percentage, the lower the variability of the data set and the inverse remains; thus, the higher the percentage obtained, the greater the variability of the data set.

It is important to note that the DE method was used in its canonical version using the DE/rand/1/bin scheme. Obviously, the comparison between multipopulation methods and DE is not strictly ideal, since DE is not a multimodal algorithm; however, this comparison is made due to the fact that DE is the base method for the formulation of MDE-ITMF and DEwI.

In stochastic optimization algorithms, an important task to obtain good results is the adequate choice of the values of the control parameters. In the case of stochastic methods devoted to obtain multiple points of global optimum simultaneously, the ideal choice of parameters is understood to be the choice that leads to obtaining the complete set of solutions in a single run, with the desired accuracy and in the minimum time possible. We will employ a sensitivity analysis method to tune the parameters for DE, MDE-ITMF and DEwI methods. This methodology fixes the value of the parameters and generates results for different values of the parameter being analyzed.

We will use the same parameters for the three methods, in order to generate an impartial environment for the simulations and results presented in this work.

As the initial step of the sensitivity analysis, it is necessary to establish the initial choice for the values of each parameter. \cite{bib:storn1997differential} suggest values for parameters $F = 0.5$ and $\emph{CR} = 0.1$ and an interval between $5 \times d$ and $10 \times d$ for the number of individuals, where $d$ is the dimension of the problem. The values of the other pertinent parameters to the MDE-ITMF method were chosen empirically from the previous knowledge about the method: $ \beta = 2\times10^{3} $, $ \rho = 1.5 $ and $ \epsilon = 5 \times 10^{-5} $. 

Table \ref{tab:variacao_NP_AS} presents the sensitivity analysis of the MDE-ITMF method for the Himmelblau function, with respect to $N_\text{P}$ (the number of individual for each subpopulation). An analysis of Table \ref{tab:variacao_NP_AS} indicates that the most suitable value to adopt is $N_{\text{p}} = 30$ (high mean value for NGP). 

\begin{table}[!htp]
\caption{Evaluation of $N_{\text{p}}$ effect in the performance of MDE-ITMF}
\label{tab:variacao_NP_AS}
\centering{}
\resizebox{\textwidth}{!}{
\begin{tabular}{cccccccccc}        
\hline
\multirow{2}{*}{$ N_{\text{p}} $} & \multicolumn{3}{c}{mean} & \multicolumn{3}{c}{standard deviation} & \multicolumn{3}{c}{coefficient of variation} \\ \cline{2-10} 
 & ET & NFE & NGP & ET & NFE & NGP & ET & NFE & NGP \\ \hline
  8 & 2,6182 & 9274,2667 & 3,3333 & 1,6019 & 9225,2711 & 0,7581 & 61,18 & 99,47 & 22,74 \\
 10 & 2,6100 & 8327,5333 & 3,7000 & 0,4028 & 2462,5795 & 0,6513 & 15,43 & 29,57 & 17,60 \\
 12 & 3,1561 & 10659,7333 & 3,7000 & 1,0651 & 6752,3554 & 0,4661 & 33,75 & 63,34 & 12,60 \\
 15 & 3,7103 & 12015,6000 & 3,6333 & 0,5804 & 3821,6052 & 0,5561 & 15,64 & 31,81 & 15,30 \\
 20 & 5,3599 & 19642,8667 & 3,8333 & 1,8523 & 12206,4054 & 0,4611 & 34,56 & 62,14 & 12,03 \\
 25 & 6,6812 & 24553,4333 & 3,8333 & 2,2476 & 15171,5929 & 0,4611 & 33,64 & 61,79 & 12,03 \\
 30 & 7,5118 & 27372,8000 & 3,9333 & 1,0680 & 6854,6230 & 0,2537 & 14,22 & 25,04 & 6,45 \\
 35 & 8,7987 & 31511,0333 & 3,8667 & 1,3218 & 8870,0095 & 0,3457 & 15,02 & 28,15 & 8,94 \\
 40 & 10,4799 & 39861,0667 & 3,9000 & 3,4938 & 23915,1964 & 0,3051 & 33,34 & 60,00 & 7,82 \\ \hline
\end{tabular} 
}
\end{table}

This sensitivity analysis was similarly applied to the parameters $F, \emph{CR}$ and $\rho$ in all functions used in this article. 


In order to determine a suitable value for the parameter $tol$ (used in DEwI), a sensitivity analysis was performed using the values $5 \times 10^{-1}, 4 \times 10^{-1}, 3 \times 10^{-1}, 2 \times 10^{-1},  10^{-1}, 10^{-2}, 10^{-3}, 5 \times 10^{-4}, 2.5 \times 10^{-4}, 10^{-4}$. The simulations of the sensitivity analysis indicate that $tol = 5 \times 10^{-4}$ presented the best set of results, considering the values for ET, NFE and NGP and the statistical measures.

Table \ref{tab:tab_parametros} presents the values of the parameters of the DE, MDE-ITMF and DEwI methods obtained with the application of sensitivity analysis in all benchmark problems analyzed here.

\begin{table}[!htp]
\centering{}
\caption{Control parameters for ED, MDE-ITMF and DEwI methods}
\resizebox{\textwidth}{!}{
\begin{tabular}{cccccccccccc} 
\hline
\multirow{3}{*}{{\small Functions}} & \multicolumn{11}{c}{Parameters} \\ \cline{2-12}
 & \multicolumn{5}{c}{{\small DE,  MDE-ITMF and DEwI}} & & \multicolumn{3}{c}{{\footnotesize MDE-ITMF and DEwI}} & & DEwI \\ \cline{2-6} \cline{8-10} \cline{12-12}
& $ G_{\text{max}} $ & $ N_{\text{p}} $ & $ F $ & $ \text{\emph{CR}} $ & $ \epsilon $ & & $ N_{\text{sp}} $ & $ \beta $ & $ \rho $ & & $ tol $ \\ \hline
$ B_{1} $ & $ 10^{3} $ & 30 & 0,7 & 0,8 & $ 5 \times 10^{-5} $ & & 4 & $ 2 \times 10^{3} $ & 2 & & $ 5 \times 10^{-4} $ \\
$ B_{2} $ & $ 10^{3} $ & 15 & 0,4 & 0,3 & $ 5 \times 10^{-5} $ & & 2 & $ 2 \times 10^{3} $ & 1 & & $ 5 \times 10^{-4} $ \\
$ B_{3} $ & $ 10^{3} $ & 20 & 0,7 & 0,8 & $ 5 \times 10^{-5} $ & & 2 & $ 2 \times 10^{3} $ & 0,6 & &  $ 5 \times 10^{-4} $ \\
$ B_{4} $ & $ 10^{3} $ & 15 & 0,6 & 0,7 & $ 5 \times 10^{-5} $ & & 4 & $ 2 \times 10^{3} $ & 0,8 & &  $ 5 \times 10^{-4} $ \\
$ B_{5} $ & $ 10^{3} $ & 30 & 0,8 & 0,7 & $ 5 \times 10^{-5} $ & & 2 & $ 2 \times 10^{3} $ & 3,2 & &  $ 5 \times 10^{-4} $ \\
$ B_{6} $ & $ 10^{3} $ & 25 & 0,6 & 0,6 & $ 5 \times 10^{-5} $ & & 3 & $ 2 \times 10^{3} $ & 2 & &  $ 5 \times 10^{-4} $ \\
$ B_{7} $ & $ 10^{3} $ & 30 & 0,6 & 0,8 & $ 5 \times 10^{-5} $ & & 4 & $ 2 \times 10^{3} $ & 0,7 & &  $ 5 \times 10^{-4} $ \\
$ B_{8} $ & $ 10^{3} $ & 20 & 0,5 & 0,3 & $ 5 \times 10^{-5} $ & & 2 &  $ 2 \times 10^{3} $ & 1,1 & &  $ 5 \times 10^{-4} $ \\
$ B_{9} $ & $ 10^{3} $ & 20 & 0,4 & 0,7 & $ 5 \times 10^{-5} $ & & 2 & $ 2 \times 10^{3} $ & 0,15 & &  $ 5 \times 10^{-4} $ \\
$ B_{10} $ & $ 10^{3} $ & 20 & 0,4 & 0,4 & $ 5 \times 10^{-5} $ & & 2 & $ 2 \times 10^{3} $ & 1,1 & &  $ 5 \times 10^{-4} $ \\ \hline
\end{tabular}%
\label{tab:tab_parametros}
}
\end{table}

All results presented were obtained using the Scilab 6.0.1 software on a computer with the following configurations: Intel(R) Core(TM) i5-4210U CPU @ 1.70GHz with 8,00 GB of RAM memory. 

Table~\ref{tab:result_media_metodos_final} presents the values of the arithmetic means for the quantities ET, NFE and NGP using the methods DE, MDE-ITMF and DEwI in the benchmark functions. It can be noted that the canonical DE algorithm shows the lowest values for $\mu_{ _{\text{ET}}}$ in the majority of cases. This fact should come as no surprise, since the canonical DE was used without the addition of any methodology to obtain different points of global minimum. Thus, with a more simple algorithm, we expect lower elapsed times. DEwI method obtained the best values for $\mu_{ _{\text{NFE}}}$ in $50\%$ of the functions. The DE and MDE-ITMF perform better in, respectively, $20\%$ and $30\%$. We can also observe a notable similarity of the behaviors for the MDE-ITMF and DEwI methods with respect to $\mu_{ _{\text{NFE}}}$. Furthermore, regarding the number of distinct global minimum points obtained, $\mu_{ _{\text{NGP}}}$, the MDE-ITMF and DEwI methods obtained all the global minimum points in, respectively, $60\%$ and $80\%$ of the applied benchmark functions. Surprisingly, the DE method obtained all points of global minimum in function B9, since DE does not use procedures that induce it to obtain distinct points of global minimum in its sequential executions.

\begin{table}[!htp]
\centering{}
\caption{Arithmetic means -- Comparative results between DE, MDE-ITMF and DEwI methods}
\resizebox{1\textwidth}{!}{%
\begin{tabular}{cccccccccccc}
\hline
\multirow{3}{*}{Functions} & \multicolumn{11}{c}{Methods} \\ \cline{2-12} 
 & \multicolumn{3}{c}{DE} & & \multicolumn{3}{c}{MDE-ITMF} & & \multicolumn{3}{c}{DEwI} \\ \cline{2-4} \cline{6-8} \cline{10-12} 
& $ \mu_{ _{\text{ET}}} $ & $ \mu_{ _{\text{NFE}}} $ & $ \mu_{ _{\text{NGP}}} $ & & $ \mu_{ _{\text{ET}}} $ & $ \mu_{ _{\text{NFE}}} $ & $ \mu_{ _{\text{NGP}}} $ & & $ \mu_{ _{\text{ET}}} $ & $ \mu_{ _{\text{NFE}}} $ & $ \mu_{ _{\text{NGP}}} $  \\ \hline 
 $ B_{1} $ & 4,3992 & 33568,50 & 1,77 &  & 6,6433 & 19315,22 &  4,00 &  & 6,6573 & 19259,5600 & 4,00 \\ 
$ B_{2} $ & 3,4202 & 40654,65 & 1,55 &  & 7,0698 & 45685,40 & 2,00 &  & 7,6328 & 46279,3800 & 2,00 \\ 
$ B_{3} $ & 1,6540 & 11646,40 & 1,45 &  & 2,0890 & 6569,48 & 2,00  &  &  2,2147 & 6631,2200 & 2,00 \\ 
$ B_{4} $ & 2,7688 & 17705,40 & 2,72 &  & 4,4484 & 10678,09 & 3,98 &  & 3,9081 & 10680,3000 & 4,00 \\ 
$ B_{5} $ & 5,4211 & 53448,90 & 1,04 &  & 3,3024 & 10858,00 & 1,96  &  & 3,4490 & 10843,3000 & 2,00 \\ 
$ B_{6} $ & 3,0739 & 21739,50 & 1,99 &  & 4,1775 & 12932,55 & 2,98 &  & 4,3893 & 12839,2700 & 2,99 \\ 
$ B_{7} $ & 5,7633 & 40290,60 & 2,80 &  & 6,4184 & 17315,76 & 4,00  &  & 6,4961 & 17324,9400 & 4,00 \\ 
$ B_{8} $ & 1,9610 & 17137,60 & 1,35 &  & 3,4193 & 16622,12 & 1,91  &  & 3,3191 & 16411,1600 & 1,98 \\ 
$ B_{9} $ & 1,4929 & 10211,80 & 2,00 &  & 2,5221 & 10557,60 & 2,00 &  & 2,5548 & 10288,4600 & 2,00 \\ 
$ B_{10} $ & 2,2514 & 18192,40 & 1,52 &  & 2,1616 & 7236,46 & 2,00 &  & 2,2718 & 7223,0600 & 2,00 \\ \hline
\label{tab:result_media_metodos_final}
\end{tabular}%
}
\end{table}

Finally, from the results presented in Table~\ref{tab:result_media_metodos_final}, it may be concluded that the DEwI method showed best values with respect to $\mu_{ _{\text{NGP}}}$ in $100\%$ of the benchmark functions, but the results of MDE-ITMF were very similar. As pointed previously, the comparison with the canonical DE was presented only to indicate the difference between a ordinary and a multimodal optimization algorithm. In short, both methods MDE-ITMF and DEwI proved to be able to locate multiple minima in the studied benchmark problems.

Table~\ref{tab:result_desvPadrao_metodos_final} presents the standard deviation for the variables ET, NFE and NGP in the benchmark functions. The standard deviation values inform how dispersed the results are in relation to the mean.

\begin{table}[!htp]
\centering{}
\caption{Standard deviations -- Comparative results between DE, MDE-ITMF and DEwI methods}
\resizebox{1\textwidth}{!}{%
\begin{tabular}{cccccccccccc}
\hline
\multirow{3}{*}{Functions} & \multicolumn{11}{c}{Methods} \\ \cline{2-12} 
 & \multicolumn{3}{c}{DE} & & \multicolumn{3}{c}{MDE-ITMF} & & \multicolumn{3}{c}{DEwI} \\ \cline{2-4} \cline{6-8} \cline{10-12} 
& $ \sigma_{ _{\text{ET}}} $ & $ \sigma_{ _{\text{NFE}}} $ & $ \sigma_{ _{\text{NGP}}} $ & & $ \sigma_{ _{\text{ET}}} $ & $ \sigma_{ _{\text{NFE}}} $ & $ \sigma_{ _{\text{NGP}}} $ & & $ \sigma_{ _{\text{ET}}} $ & $ \sigma_{ _{\text{NFE}}} $ & $ \sigma_{ _{\text{NGP}}} $  \\ \hline
$ B_{1} $ & 0,4592 & 7035,5001 & 0,5835 &  & 0,2576 & 749,6291 &  0,0000 &  &  0,1882 &  669,1326 &  0,0000 \\ 
$ B_{2} $  & 1,8940 & 27412,8026 & 0,5000 &  & 0,7977 & 5453,6575 &  0,0000 &  &  0,6450 &  4255,5780 &  0,0000 \\ 
$ B_{3} $  & 0,8993 & 12249,1390 & 0,5000 &  & 0,1562 &  324,7391 &  0,0000 &  &  0,1473 & 362,6995 &  0,0000 \\ 
$ B_{4} $ & 0,9777 & 11414,6090 & 0,6526 &  & 0,7930 & 1109,7794 & 0,2000 &  &  0,1736 &  534,7423 &  0,0000 \\ 
$ B_{5} $ & 3,2716 & 41167,2972 & 0,5490 &  &  0,2130 & 627,4527 & 0,1969 &  &  0,1443 &  540,6486 &  0,0000 \\ 
$ B_{6} $ & 0,9867 & 13907,8584 & 0,6113 &  & 0,3574 &  1872,9067 & 0,2000 &  &  0,3164 & 2027,0203 &  0,1000 \\ 
$ B_{7} $ & 1,1085 & 13550,9537 & 0,6816 &  &  0,2559 & 614,7779 &  0,0000 &  & 0,3309 &  596,0172 &  0,0000 \\ 
$ B_{8} $ &  0,1468 &  1661,4325 & 0,4794 &  & 0,3226 & 2224,8091 & 0,2876 &  & 0,2757 & 2029,1889 &  0,1407 \\ 
$ B_{9} $ &  0,1012 &  827,9857 &  0,0000 &  & 0,2170 & 1033,0543 &  0,0000 &  & 0,1741 & 897,8128 &  0,0000 \\ 
$ B_{10} $ & 1,7495 & 22250,4364 & 0,5021 &  & 0,1620 & 560,7570 &  0,0000 &  &  0,1482 &  448,3346 &  0,0000 \\ \hline
\label{tab:result_desvPadrao_metodos_final}
\end{tabular}%
}
\end{table}

The DE method obtained the lowest values for $\sigma_{ _{\text{ET}}}$ and $\sigma_{ _{\text{NFE}}}$ in $20\%$ of the functions. This method presented the smallest $\sigma_{ _{\text{NGP}}}$ in function $B_{9}$ in agreement with the value presented in Table~\ref{tab:result_media_metodos_final}. The MDE-ITMF method presented the lowest values of $\sigma_{ _{\text{ET}}}$ in $10\%$ of the analyzed functions and $20\%$ in relation to $\sigma_{ _{\text{NFE}}}$. The superiority of MDE-ITMF over the canonical DE (as a multimodal algorithm) can be evaluated by $\sigma_{ _{\text{NGP}}}$. This result is consistent to that presented in Table~\ref{tab:result_media_metodos_final}. The standard deviation values were equal to zero in the functions where the MDE-ITMF method obtained all the points of global minimum. The DEwI method presented the lowest values for the standard deviations in comparison with the canonical DE and MDE-ITMF. When evaluating the values of $\sigma_{ _{\text{ET}}}$ it is noted that the method was successful compared to the other methods in $70\%$ of the functions. Regarding $\sigma_{ _{\text{NFE}}}$, success occurred in $60\%$ of the functions. The real gain of the DEwI method can be analyzed with the values for NGP; in this case the method obtained the best values of $\sigma_{ _{\text{NGP}}}$ in $100\%$ of the functions, with $\sigma_{ _{\text{NGP}}}$ equal to zero in eight of ten benchmark functions.

Table~\ref{tab:result_coefVaria_metodos_final} presents the coefficient of variation for the variables ET, NFE and NGP. The results presented in Table~\ref{tab:result_coefVaria_metodos_final} are in accordance with the aforementioned values presented in the Table~\ref{tab:result_desvPadrao_metodos_final}. This was already expected taking into account that the coefficient of variation is the ratio between the values of the standard deviation and the values of the arithmetic mean. In this way, when the standard deviation has a very small value, this causes the coefficient of variation to present a low value, since the arithmetic mean is a non-zero number. The coefficient of variation measures the variability of the data set, by obtaining low values for the coefficient of variation it can be said that the data sets of the solutions are homogeneous. Once again, we can note that DEwI proved superior to MDE-ITMF, but with a little diference.

\begin{table}[!htp]
\centering{}
\caption{Coefficient of variation -- Comparative results between DE, MDE-ITMF and DEwI methods}
\resizebox{\textwidth}{!}{%
\begin{tabular}{cccccccccccc}
\hline
\multirow{3}{*}{Functions} & \multicolumn{11}{c}{Methods} \\ \cline{2-12} 
 & \multicolumn{3}{c}{DE} & & \multicolumn{3}{c}{MDE-ITMF} & & \multicolumn{3}{c}{DEwI} \\ \cline{2-4} \cline{6-8} \cline{10-12} 
 & $ \upsilon_{ _{\text{ET}}} (\%) $ & $ \upsilon_{ _{\text{NFE}}} (\%)$ & $ \upsilon_{ _{\text{NGP}}} (\%) $ &  & $ \upsilon_{ _{\text{ET}}} (\%) $ & $ \upsilon_{ _{\text{NFE}}} (\%)$ & $ \upsilon_{ _{\text{NGP}}} (\%) $ &  & $ \upsilon_{ _{\text{ET}}} (\%) $ & $ \upsilon_{ _{\text{NFE}}} (\%) $ & $ \upsilon_{ _{\text{NGP}}} (\%) $ \\ \hline
$ B_{1} $ & 10,44 & 20,96 & 32,97 &  & 3,88 & 3,88 &  0,00 &  &  2,83 &  3,47 &  0,00 \\ 
$ B_{2} $ & 55,38 & 67,43 & 32,26 &  & 11,28 & 11,94 &  0,00 &  &  8,45 &  9,20 &  0,00 \\ 
$ B_{3} $ & 54,37 & 105,18 & 34,48 &  & 7,48 &  4,94 &  0,00 &  &  6,65 & 5,47 &  0,00 \\ 
$ B_{4} $ & 35,31 & 64,47 & 23,99 &  & 17,83 & 10,39 & 5,03 &  &  4,44 &  5,01 &  0,00 \\ 
$ B_{5} $ & 60,35 & 77,02 & 52,79 &  & 6,45 & 5,78 & 10,05 &  &  4,18 &  4,99 &  0,00 \\ 
$ B_{6} $ & 32,10 & 63,98 & 30,72 &  & 8,56 &  14,48 & 6,71 &  &  7,21 & 15,79 &  3,34 \\ 
$ B_{7} $ & 19,23 & 33,63 & 24,34 &  &  3,99 & 3,55 &  0,00 &  & 5,09 &  3,44 &  0,00 \\ 
$ B_{8} $ &  7,49 &  9,69 & 35,51 &  & 9,43 & 13,38 & 15,06 &  & 8,31 & 12,36 &  7,11 \\ 
$ B_{9} $ &  6,78 &  8,11 &  0,00 &  & 8,60 & 9,78 &  0,00 &  & 6,82 & 8,73 &  0,00 \\ 
$ B_{10} $ & 77,71 & 122,31 & 33,03 &  & 7,50 & 7,75 &  0,00 &  &  6,53 &  6,21 &  0,00 \\ \hline
\end{tabular}%
}
\label{tab:result_coefVaria_metodos_final}
\end{table}



\section{Conclusions} \label{sec:conclusions}

In this work we proposed two new multimodal and multipopulation methods based on the well-known Differential Evolution method, called Multipopulation Differential Evolution with Iterative Technique of Modification of the Objective Function (MDE-ITMF) and Differential Evolution with Initialization (DEwI). Both methods were tested in a set of benchmark functions and the results were compared with respect of elapsed times, number of function evaluations and number of global points found.


The MDE-ITMF method was able to be applied to multimodal problems, obtaining the best values for the arithmetic mean of the number of distinct global minimum points in $60\%$ of the benchmark functions. However, in the other functions, the method also obtained satisfactory values for the number of global minima found. The DEwI method has shown the best performance with respect to the number of distinct global optimum points obtained in a single run.

The computational results indicated that both techniques are capable to obtain, in a reliable way and with low computational cost, the full set of minima in the benchmark functions.

Furthermore, for both MDE-ITMF and DEwI, the elapsed time can be reduced parallelizing the methods with respect to the generation of the initial populations and the applications of the Differential Evolution procedures. Thus, both methods look promising in future applications, such as real-world problems.


%
%

\begin{acknowledgements}
This study was financed in part by the Coordenação de A\-{per}\-{fei}\-{ço}\-{a}\-{men}\-{to} de Pessoal de Nível Superior - Brasil (CAPES) - Finance Code 001.
\end{acknowledgements}

%
%

\bibliographystyle{spbasic}      
\bibliography{references}   

%
%

\end{document}